\documentclass[reqno,10pt]{amsart}
\usepackage{amsmath,amssymb}
\usepackage{latexsym}
\usepackage[dvips]{graphics}
\usepackage{epsfig}

\newtheorem{thm}{Theorem}[section]
\newtheorem{lem}[thm]{Lemma}

\theoremstyle{definition}
\newtheorem{rek}[thm]{Remark}

\newcommand{\ncr}[2]{{#1 \choose #2}}

\newcommand\ben{\begin{enumerate}}
\newcommand\een{\end{enumerate}}

\newcommand{\C}{\mathbb{C}}

\newcommand{\ga}{\alpha}    
\newcommand{\gb}{\beta}      
\newcommand{\g}{\gamma}      

\newcommand{\twocase}[5]{#1 \begin{cases} #2 & \text{\rm #3}\\ #4
&\text{\rm #5} \end{cases}  }

\newcommand{\threecase}[7]{#1 \begin{cases} #2 & \text{\rm #3}\\ #4
&\text{\rm #5}\\ #6 & \text{\rm #7} \end{cases}  }

\newcommand\be{\begin{equation}}
\newcommand\ee{\end{equation}}
\newcommand\bea{\begin{eqnarray}}
\newcommand\eea{\end{eqnarray}}

\newcommand{\foh}{\frac{1}{2}}  

\renewcommand{\d}{{\mathrm{d}}} 

\newcommand{\rs}{{\rm RS}}
\newcommand{\ra}{{\rm RA}}
\newcommand{\rso}{{\rm RS_{\rm obs}}}
\newcommand{\rao}{{\rm RA_{\rm obs}}}

\numberwithin{equation}{section}

\begin{document}

\title{A Derivation of the Pythagorean Won-Loss Formula in Baseball}

\author{Steven J. Miller}
\address{Department of Mathematics, Brown University, 151 Thayer
 Street,
Providence, RI 02912}
 \email{sjmiller@math.brown.edu}



\date{\today}

\begin{abstract}
It has been noted that in many professional sports leagues a good
predictor of a team's end of season won-loss percentage is Bill
James' Pythagorean Formula $\frac{\rso^\g}{\rso^\g+\rao^\g}$,
where $\rso$ (resp. $\rao$) is the observed average number of runs
scored (allowed) per game and $\gamma$ is a constant for the
league; for baseball the best agreement is when $\gamma$ is about
$1.82$. This formula is often used in the middle of a season to
determine if a team is performing above or below expectations, and
estimate their future standings.

We provide a theoretical justification for this formula and value of
$\gamma$ by modeling the number of runs scored and allowed in
baseball games as independent random variables drawn from Weibull
distributions with the same $\beta$ and $\gamma$ but different
$\ga$; the probability density is \be \twocase{f(x;\ga,\gb,\g) \ = \
}{\frac{\g}{\ga}\ ((x-\gb)/\ga)^{\g-1}\ e^{- ((x-\gb)/\ga)^{\g}}}{if
$x \ge \gb$}{0}{otherwise.} \nonumber\ \ee This model leads to a
predicted won-loss percentage of $\frac{(\rs-\gb)^\g}{(\rs-\gb)^\g +
(\ra-\gb)^\g}$; here $\rs$ (resp. $\ra$) is the mean of the Weibull
random variable corresponding to runs scored (allowed), and $\rs -
\beta$ (resp. $\ra - \beta$) is an estimator of $\rso$ (resp.
$\rao$). An analysis of the 14 American League teams from the 2004
baseball season shows that (1) given that the runs scored and
allowed in a game cannot be equal, the runs scored and allowed are
statistically independent; (2) the best fit Weibull parameters
attained from a least squares analysis and the method of maximum
likelihood give good fits. Specifically, least squares yields a mean
value of $\g$ of $1.79$ (with a standard deviation of $.09$) and
maximum likelihood yields a mean value of $\g$ of $1.74$ (with a
standard deviation of $.06$), which agree beautifully with the
observed best value of $1.82$ attained by fitting
$\frac{\rso^\g}{\rso^\g+\rao^\g}$ to the observed winning
percentages.
\end{abstract}

\maketitle


\section{Introduction}

The goal of this paper is to derive Bill James' Pythagorean Formula
(see \cite{James}, as well as \cite{Angus,Oliver}) from reasonable
assumptions about the distribution of scores. Given a sports league,
if the observed average number of runs a team scores and allows are
$\rso$ and $\rao$, then the Pythagorean Formula predicts the team's
won-loss percentage should be $\frac{\rso^\g}{\rso^\g+\rao^\g}$ for
some $\g$ which is constant for the league. Initially in baseball
the exponent $\g$ was taken to be $2$ (which led to the name),
though fitting $\g$ to the observed records from many seasons lead
to the best $\g$ being about $1.82$. Often this formula is applied
part way through a season to estimate a team's end of season
standings. For example, if halfway through a season a team has far
more wins than this formula predicts, analysts often claim the team
is playing over their heads and predict they will have a worse
second-half.

Rather than trying to find the best $\gamma$ by looking at many
teams' won-loss percentages, we take a different approach and derive
the formula and optimal value of $\gamma$ by modeling the runs
scored and allowed each game for a team as independent random
variables drawn from Weibull distributions with the same $\gb$ and
$\g$ but different $\ga$ (see \S\ref{sec:AL2004} for an analysis of
the 2004 season which shows that, subject to the condition that the
runs scored and allowed in a game must be distinct integers, the
runs scored and allowed are statistically independent, and
\S\ref{sec:concfutwork} for additional comments on the
independence). Recall the three-parameter Weibull distribution (see
also \cite{Fe2}) is \be \twocase{f(x;\ga,\gb,\g) \ = \
}{\frac{\g}{\ga} \left(\frac{x-\gb}{\ga}\right)^{\g-1} e^{-
((x-\gb)/\ga)^{\g}}}{if $x \ge \gb$}{0}{otherwise.} \ee We denote
the means by $\rs$ and $\ra$, and we show below that $\rs - \beta$
(resp. $\ra - \beta$) is an estimator of the observed average number
of runs scored (resp. allowed) per game. The reason $\rs - \beta$
and not $\rs$ is the estimator of the observed average runs scored
per game is due to the discreteness of the runs scored data; this is
described in greater detail below. Our main theoretical result is
proving that this model leads to a predicted won-loss percentage of
\be\label{eq:wonlossperc} \mbox{\rm Won-Loss
Percentage}(\rs,\ra,\gb,\g) \ = \ \frac{(\rs-\gb)^\g}{(\rs-\gb)^\g +
(\ra-\gb)^\g}; \ee note for all $\g$ that if $\rs = \ra$ in
\eqref{eq:pythagammap} then as we would expect the won-loss
percentage is $50\%$.

In \S\ref{sec:AL2004} we analyze in great detail the 2004 baseball
season for the 14 teams of the American League. Complete results of
each game are readily available (see for example \cite{Almanac}),
which greatly facilitates curve fitting and error analysis. For each
of these teams we used the method of least squares and the method of
maximum likelihood to find the best fit Weibulls to the runs scored
and allowed per game (with each having the same $\gamma$ and both
having $\beta = -.5$; we explain why this is the right choice for
$\beta$ below). Standard $\chi^2$ tests (see for example
\cite{CaBe}) show our fits are adequate. For continuous random
variables representing runs scored and runs allowed, there is zero
probability of both having the same value; the situation is markedly
different in the discrete case. In a baseball game runs scored and
allowed \emph{cannot} be entirely independent, as games do not end
in ties; however, modulo this condition, modified $\chi^2$ tests
(see \cite{BF,SD}) do show that, given that runs scored and allowed
per game must be distinct integers, the runs scored and allowed per
game are statistically independent. See \cite{Ci} for more on the
independence of runs scored and allowed.

Thus the assumptions of our theoretical model are met, and the
Pythagorean Formula should hold for some exponent $\g$. Our main
experimental result is that, averaging over the 14 teams, the method
of least squares yields a mean of $\gamma$ of $1.79$ with a standard
deviation of $.09$ (the median is $1.79$ as well); the method of
maximum likelihood yields a mean of $\gamma$ of $1.74$ with a
standard deviation of $.06$ (the median is $1.76$). This is in line
with the numerical observation that $\gamma = 1.82$ is the best
exponent.

In order to obtain simple closed form expressions for the
probability of scoring more runs than allowing in a game, we assume
that the runs scored and allowed are drawn from continuous and not
discrete distributions. This allows us to replace discrete sums with
continuous integrals, and in general integration leads to more
tractable calculations than summations. Of course assumptions of
continuous run distribution cannot be correct in baseball, but the
hope is that such a computationally useful assumption is a
reasonable approximation to reality; it may be more reasonable in a
sport such as basketball, and this would make an additional,
interesting project. Closed form expressions for the mean, variance
and probability that one random variable exceeds another are
difficult for general probability distributions; however, the
integrations that arise from a Weibull distribution with parameters
$(\ga,\gb,\g)$ are very tractable. Further, as the three parameter
Weibull is a very flexible family and takes on a variety of
different shapes, it is not surprising that for an appropriate
choice of parameters it is a good fit to the runs scored (or
allowed) per game. What is fortunate is that we can get good fits to
both runs scored and allowed simultaneously, using the same $\g$ for
each; see \cite{BFAM} for additional problems modeled with Weibull
distributions. For example, $\g=1$ is the exponential and $\g = 2$
is the Rayleigh distribution. Note the great difference in behavior
between these two distributions. The exponential's maximum
probability is at $x=\gb$, whereas the Rayleigh is zero at $x=\gb$.
Additionally, for any $M > \gb$ any Weibull has a non-zero
probability of a team scoring (or allowing) more than $M$ runs,
which is absurd of course in the real world. The tail probabilities
of the exponential are significantly greater than those of the
Rayleigh, which indicates that perhaps something closer to the
Rayleigh than the exponential is the truth for the distribution of
runs.

We have incorporated a translation parameter $\gb$ for several
reasons. First, to facilitate applying this model to sports other
than baseball. For example, in basketball no team scores fewer than
20 points in a game, and it is not unreasonable to look at the
distribution of scores above a baseline. A second consequence of
$\gb$ is that adding $P$ points to both the runs scored and runs
allowed each game does not change the won-loss percentage; this is
reflected beautifully in \eqref{eq:wonlossperc}, and indicates that
it is more natural to measure scores above a baseline (which may be
zero). Finally, and most importantly, as remarked there are issues
in the discreteness of the data and the continuity of the model. In
the least squares and maximum likelihood curve fitting we bin the
runs scored and allowed data into bins of length $1$; for example, a
natural choice of bins is \be\label{eq:binsforWMLSdontuse} [0, 1)\
\cup\ [1, 2)\ \cup\ \cdots\ \cup\ [9, 10)\ \cup\ [10, 12) \ \cup \
[12, \infty). \ee As baseball scores are non-negative integers, all
of the mass in each bin is at the left endpoint. If we use
untranslated Weibulls (i.e., $\gb=0$) there would be a discrepancy
in matching up the means.

For example, consider a simple case when in half the games the team
scores 0 runs and in the other half they score 1. Let us take as our
bins $[0,1)$ and $[1,2)$, and for ease of exposition we shall find
the best fit function constant on each bin. Obviously we take our
function to be identically $\foh$ on $[0,2)$; however, the observed
mean is $\foh \cdot 0 + \foh \cdot 1 = \foh$ whereas the mean of our
piecewise constant approximant is $1$. If instead we chose
$[-.5,.5)$ and $[.5,1.5)$ as our bins then the approximant would
also have a mean of $\foh$. Returning to our model, we see a better
choice of bins is \be\label{eq:listbaseballbins} [-.5, .5] \ \cup \
[.5, 1.5] \ \cup \ \cdots \ \cup \ [7.5, 8.5]\ \cup\ [8.5, 9.5]\
\cup\ [9.5, 11.5] \ \cup\ [11.5, \infty). \ee An additional
advantage of the bins of \eqref{eq:listbaseballbins} is that we may
consider either open or closed endpoints, as there are no baseball
scores that are half-integral. Thus, in order to have the baseball
scores in the \emph{center} of their bins, we take $\gb = -.5$ and
use the bins in \eqref{eq:listbaseballbins}. In particular, if the
mean of the Weibull approximating the runs scored (resp. allowed)
per game is $\rs$ (resp. $\ra$) then $\rs-\gb$ (resp. $\ra - \gb$)
is an estimator of the observed average number of runs scored (resp.
allowed) per game.


\section{Theoretical Model and Predictions}\label{sec:theorymodel}

We determine the mean of a Weibull distribution with parameters
$(\ga,\gb,\g)$, and then use this to prove our main result, the
Pythagorean Formula (Theorem \ref{thm:pythmweibull}). Let
$f(x;\ga,\gb,\g)$ be the probability density of a Weibull with
parameters $(\ga,\gb,\g)$: \be \twocase{f(x;\ga,\gb,\g) \ = \
}{\frac{\g}{\ga} \left(\frac{x-\gb}{\ga}\right)^{\g-1}
e^{-((x-\gb)/\ga)^\g}}{if $x\ge \gb$}{0}{otherwise.} \ee For $s\in
\C$ with the real part of $s$ greater than $0$, recall the
$\Gamma$-function (see \cite{Fe1}) is defined by
\begin{equation}
\Gamma(s)\ =\ \int_0^\infty e^{-u} u^{s-1} \d u \ = \
\int_0^\infty e^{-u} u^s \frac{\d u}{u}.
\end{equation} Letting $\mu_{\ga,\gb,\g}$ denote the mean of $f(x;\ga,\gb,\g)$,
we have \bea \mu_{\ga,\gb,\g} & \ = \ & \int_\gb^\infty x \cdot
\frac{\g}{\ga} \left(\frac{x-\gb}{\ga}\right)^{\g-1}
e^{-((x-\gb)/\ga)^\g}\d x\nonumber\\ & = & \int_\gb^\infty  \ga
\frac{x-\gb}{\ga} \cdot \frac{\g}{\ga}
\left(\frac{x-\gb}{\ga}\right)^{\g-1}e^{-((x-\gb)/\ga)^\g}\d x\ +\
\gb. \eea We change variables by setting $u =
\left(\frac{x-\gb}{\ga}\right)^\g$. Then $\d u = \frac{\g}{\ga}
\left(\frac{x-\gb}{\ga}\right)^{\g-1}\d x$ and we have \bea
\mu_{\ga,\gb,\g} & \ = \ & \int_0^\infty \ga u^{\gamma^{-1}} \cdot
e^{-u} \d u \ + \ \gb \nonumber\\ & = & \ga \int_0^\infty e^{-u}
u^{1+\g^{-1}} \frac{\d u}{u} \ + \ \gb \nonumber\\ & = & \ga
\Gamma(1+\g^{-1}) \ + \ \gb.\eea

A similar calculation determines the variance. We record these
results:

\begin{lem}\label{lem:Weibullmeanvar} The mean $\mu_{\ga,\gb,\g}$ and
variance $\sigma^2_{\ga,\gb,\g}$ of a Weibull with parameters
$(\ga,\gb,\g)$ are \bea \mu_{\ga,\gb,\g} &\ =\ & \ga
\Gamma(1+\g^{-1}) + \gb \nonumber\\ \sigma^2_{\ga,\gb,\g}&\ =\ &
\ga^2 \Gamma\left(1+2\g^{-1}\right)  - \ga^2
\Gamma\left(1+\g^{-1}\right)^2.\eea
\end{lem}

We can now prove our main result:

\begin{thm}[Pythagorean Won-Loss Formula]\label{thm:pythmweibull}
Let the runs scored and runs allowed
per game be two independent random variables drawn from Weibull
distributions with parameters $(\ga_\rs,\gb,\g)$ and
$(\ga_\ra,\gb,\g)$ respectively, where $\ga_\rs$ and $\ga_\ra$ are
chosen so that the means are $\rs$ and $\ra$. If $\g  > 0$ then
\be\label{eq:pythagammap} \mbox{\rm Won-Loss
Percentage}(\rs,\ra,\gb,\g) \ = \ \frac{(\rs-\gb)^\g}{(\rs-\gb)^\g
+ (\ra-\gb)^\g}. \ee
\end{thm}

\begin{proof}

Let $X$ and $Y$ be independent random variables with Weibull
distributions $(\ga_\rs,\gb,\g)$ and $(\ga_\ra,\gb,\g)$
respectively, where $X$ is the number of runs scored and $Y$ the
number of runs allowed per game. As the means are $\rs$ and $\ra$,
by Lemma \ref{lem:Weibullmeanvar} we
have \bea \rs\ & \ = \ & \ga_\rs \Gamma(1+\g^{-1}) +\gb \nonumber\\
\ra & \ = \ & \ga_\ra \Gamma(1+\g^{-1}) + \gb. \eea Equivalently,
we have \bea\label{eq:genweibmeans} \ga_\rs & \ = \ &
\frac{\rs-\gb}{\Gamma(1+\g^{-1})} \nonumber\\
\ga_\ra & \ = \ & \frac{\ra-\gb}{\Gamma(1+\g^{-1})}.\eea

We need only calculate the probability that $X$ exceeds $Y$. Below
we constantly use the integral of a probability density is $1$. We
have \bea & & \mbox{Prob}(X
> Y) \ = \ \int_{x=\gb}^\infty \int_{y=\gb}^x
f(x;\ga_\rs,\gb,\g) f(y;\ga_\ra,\gb,\g) \d y\; \d x \nonumber\\
& & = \ \int_{x=\gb}^\infty\int_{y=\gb}^x \frac{\g}{\ga_\rs}
\left(\frac{x-\gb}{\ga_{RS}}\right)^{\g-1}
e^{-((x-\gb)/\ga_\rs)^\g}
\frac{\g}{\ga_\ra}\left(\frac{y-\gb}{\ga_{\ra}}\right)^{\g-1} e^{-
((y-\gb)/\ga_\ra)^\g} \d y\; \d x \nonumber\\ & & = \
\int_{x=0}^\infty\frac{\g}{\ga_\rs}
\left(\frac{x}{\ga_{RS}}\right)^{\g-1} e^{-(x/\ga_\rs)^\g} \left[
\int_{y=0}^{x}
\frac{\g}{\ga_\ra}\left(\frac{y}{\ga_{\ra}}\right)^{\g-1} e^{-
(y/\ga_\ra)^\g} \d y \right] \d x  \nonumber\\ & & = \
\int_{x=0}^\infty\frac{\g}{\ga_\rs}
\left(\frac{x}{\ga_{RS}}\right)^{\g-1} e^{-(x/\ga_\rs)^\g} \left[1
- e^{-(x/\ga_{\ra})^\g}\right] \d x  \nonumber\\ & & = \ 1 -
\int_{x=0}^\infty\frac{\g}{\ga_\rs}
\left(\frac{x}{\ga_{RS}}\right)^{\g-1} e^{-(x/\ga)^\g}\d x,\eea
where we have set \be \frac1{\ga^\g} \ = \ \frac1{\ga_\rs^\g} +
\frac{1}{\ga_\ra^\g} \ = \ \frac{\ga_\rs^\g +
\ga_\ra^\g}{\ga_\rs^\g \ga_\ra^\g}.\ee Therefore
\bea\label{eq:derivweibpythag1} \mbox{Prob}(X
> Y) & \ = \ & 1 - \frac{\ga^\g}{\ga_\rs^\g} \int_{0}^\infty
\frac{\g}{\ga}
\left(\frac{x}{\ga}\right)^{\g-1} e^{(x/\ga)^\g} \d x \nonumber\\
& = & 1 - \frac{\ga^\g}{\ga_\rs^\g} \nonumber\\ & = & 1 -
\frac1{\ga_\rs^\g} \frac{\ga_\rs^\g \ga_\ra^\g}{\ga_\rs^\g +
\ga_\ra^\g} \nonumber\\ & = & \frac{\ga_\rs^\g}{\ga_\rs^\g +
\ga_\ra^\g}.\eea Substituting the relations for $\ga_\rs$ and
$\ga_\ra$ of \eqref{eq:genweibmeans} into
\eqref{eq:derivweibpythag1} yields \bea \mbox{Prob}(X
> Y) & \ = \ &  \frac{(\rs-\gb)^\g}{(\rs-\gb)^\g + (\ra-\gb)^\g}, \eea
which completes the proof of Theorem \ref{thm:pythmweibull}.
\end{proof}

\begin{rek} The reason the integrations can be so easily
performed (determining the normalization constants, the mean and
variance, as well as calculating the probability that $X$ exceeds
$Y$) is that we have terms such as $e^{-u^\g} u^{\g-1}$; these are
very easy to integrate. It is essential, however, that we also
have a tractable expression for the mean in terms of the
parameters. Fortunately this is possible as the mean is a simple
combination of the $\Gamma$-function and the parameters. As we fix
$\g$ and then choose $\ga_\rs$ or $\ga_\ra$, it is important that
the argument of the $\Gamma$-function only involve $\g$ and not
$\ga_\rs$ or $\ga_\ra$. If the argument of the $\Gamma$-function
involved $\ga_\rs$ or $\ga_\ra$, then we would have to solve
equations of the form $\rs = g(\ga_\rs,\g)\Gamma(h(\ga_\rs,\g))$
for some functions $g$ and $h$. Inverting this to solve for
$\ga_\rs$ as a function of $\g$ and $\rs$ would be difficult in
general. Finally we remark that the essential aspect of
\eqref{eq:genweibmeans} is that $\ga_\rs$ is proportional to
$\rs-\gb$. It does not matter that the proportionality constant
involves $\g$. While it is difficult to solve $\Gamma(1+\g^{-1}) =
z$ for $\g$, we do not need to; these factors cancel.
\end{rek}

\begin{rek}
We take $\g > 0$ as if $\g < 0$ then \eqref{eq:pythagammap} (while
still true) is absurd. For example, if $\g = -.5$, $\beta = 0$, $\rs
= 25$ and $\ra = 16$, then \eqref{eq:pythagammap} predicts a winning
percentage of \be \frac{25^{-1/2}}{25^{-1/2} + 16^{-1/2}}\ =\
\frac49 \ < \ \foh; \ee thus a team that scores more runs than it
allows is predicted to have a losing season! Of course, when $\g \le
0$ we have a very strange probability distribution. Not only is the
behavior near $x=0$ interesting but we no longer have rapid decay at
infinity (the probability now falls off as $x^{\g-1}$), and this is
unlikely to be a realistic model. \end{rek}


\section{Numerical Results: American League
$2004$}\label{sec:AL2004}

We analyzed the 14 teams\footnote{The teams are ordered by division
(AL East, AL Central, AL West) and then by number of regular season
wins, with the exception of the Boston Red Sox who as World Series
champions are listed first.} of the American League from the 2004
season in order to determine the reasonableness of the assumptions
in our model; we leave the National League teams as an exercise to
the reader. We used the method of least squares\footnote{We
minimized the sum of squares of the error from the runs scored data
plus the sum of squares of the error from the runs allowed data; as
$\gb = -.5$ there were three free parameters: $\ga_\rs,\ga_\ra$ and
$\g$. Specifically, let ${\rm Bin}(k)$ be the
$k$\textsuperscript{th} bin from \eqref{eq:listbaseballbins}. If
$\rso(k)$ (resp. $\rao(k)$) denotes the observed number of games
with the number of runs scored (allowed) in ${\rm Bin}(k)$, and
$A(\ga,\gb,\g,k)$ denotes the area under the Weibull with parameters
$(\ga,\gb,\g)$ in ${\rm Bin}(k)$, then for each team we found the
values of $(\ga_\rs,\ga_\ra,\g)$ that minimized \be
\sum_{k=1}^{\#{\rm Bins}} \left( \rso(k) - \#{\rm Games} \cdot
A(\ga_\rs,-.5,\g,k)\right)^2 + \sum_{k=1}^{\#{\rm Bins}} \left(
\rao(k) - \#{\rm Games} \cdot A(\ga_\ra,-.5,\g,k)\right)^2. \ee
}\label{foot:mlsqs} and the method of maximum
likelihood\footnote{Notation as in Footnote 2, the likelihood
function of the sample is \bea L(\alpha_\rs, \alpha_\ra, -.5,
\gamma)& \ = \ & \ncr{\#{\rm Games}}{\rso(1), \dots, \rso(\#{\rm
Bins})} \prod_{k=1}^{\#{\rm Bins}} A(\ga_\rs,-.5,\gamma,k)^{\rso(k)}
\nonumber\\ & & \ \ \cdot \ \ncr{\#{\rm Games}}{\rao(1), \dots,
\rao(\#{\rm Bins})} \prod_{k=1}^{\#{\rm Bins}}
A(\ga_\ra,-.5,\gamma,k)^{\rao(k)}. \eea For each team we find the
values of the parameters $\ga_\rs$, $\ga_ra$ and $\gamma$ that
maximize the likelihood. Computationally, it is equivalent to
maximize the logarithm of the likelihood, and we may ignore the
multinomial coefficients are they are independent of the
parameters.} with the bins of \eqref{eq:listbaseballbins}. For each
team we simultaneously found the best fit Weibulls of the form
$(\ga_\rs,-.5,\g)$ and $(\ga_\ra,-.5,\g)$. We then compared the
predicted number of wins, losses, and won-loss percentage with the
actual data:

\begin{center} Results from the Method of Least Squares
\scalebox{.7}{\includegraphics{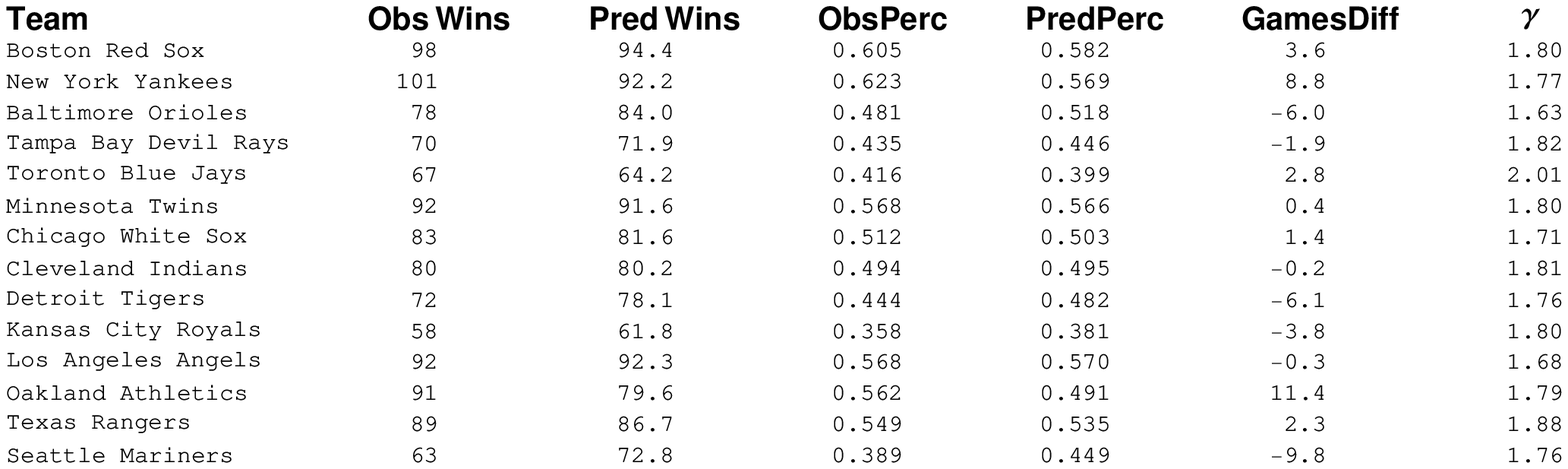}}
\end{center}

\bigskip

\begin{center} Results from the Method of Maximum Likelihood
\scalebox{.7}{\includegraphics{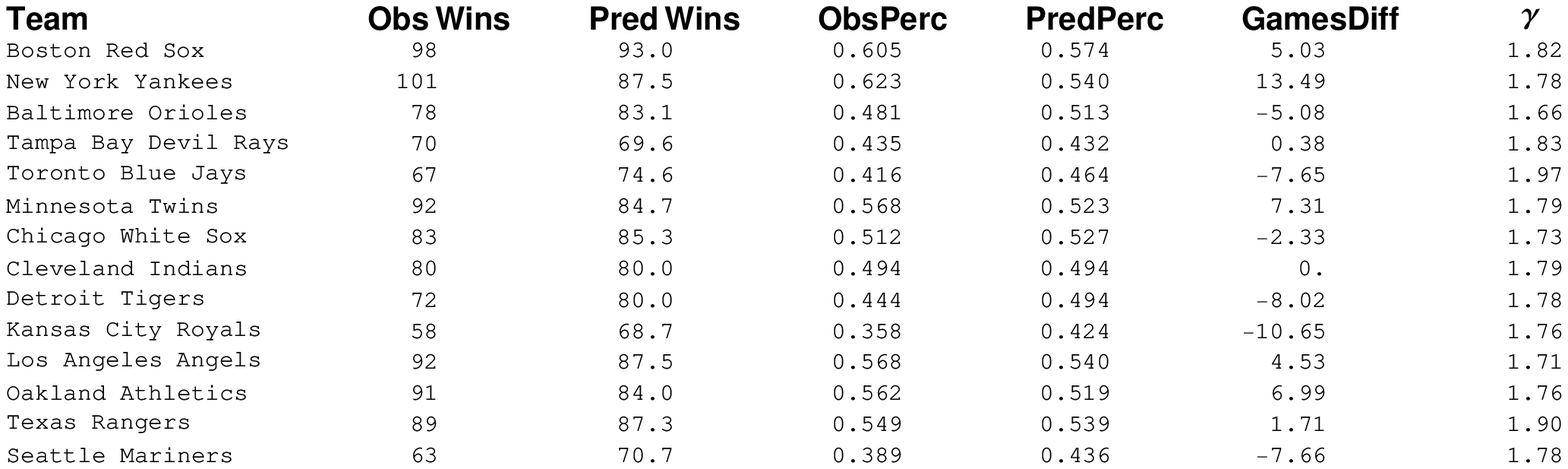}}
\end{center}

Using the method of least squares, the mean of $\g$ over the $14$
teams is $1.79$ with a standard deviation is $.09$ (the median is
$1.79$); using the method of maximum likelihood the mean of $\g$
over the $14$ teams is $1.74$ with a standard deviation of $.06$
(the median is $1.76$). Note that the numerically observed best
exponent of $1.82$ is well within this region for both approaches.

We now consider how close the estimates of team performance are to
the observed season records. For the method of least squares, over
the $14$ teams the mean number of the difference between observed
and predicted wins was $0.19$ with a standard deviation of $5.69$
(and a median of $0.07$); if we consider just the absolute value of
the difference then we have a mean of $4.19$ with a standard
deviation of $3.68$ (and a median of $3.22$). For the method of
maximum likelihood, over the $14$ teams the mean number of the
difference between observed and predicted wins was $-0.13$ with a
standard deviation of $7.11$ (and a median of $0.19$); if we
consider just the absolute value of the difference then we have a
mean of $5.77$ with a standard deviation of $3.85$ (and a median of
$6.04$). This is consistent with the observation that the
Pythagorean Formula is usually accurate to about four games in a 162
game season.

For the remainder of the paper, we analyze the fits from the method
of maximum likelihood; these fits were slightly better than those
from the method of least squares. The estimates from the method of
maximum likelihood enjoy many desirable properties, including being
asymptotically minimum variance unbiased estimators and yielding
sufficient estimators (whenever they exist).

We performed $\chi^2$ tests to determine the goodness of the fit
from the best fit Weibulls from the method of maximum
likelihood\footnote{Using the bins from \eqref{eq:listbaseballbins}
(and the rest of the notation as in Footnote 2), we studied \be
\sum_{k=1}^{\#{\rm Bins}} \frac{\left( \rso(k) - \#{\rm Games} \cdot
A(\ga_\rs,-.5,\g,k)\right)^2 }{ \#{\rm Games} \cdot
A(\ga_\rs,-.5,\g,k)} + \sum_{k=1}^{\#{\rm Bins}} \frac{\left(
\rao(k) - \#{\rm Games} \cdot A(\ga_\ra,-.5,\g,k)\right)^2 }{ \#{\rm
Games} \cdot A(\ga_\ra,-.5,\g,k)}. \ee This has a $\chi^2$
distribution with $2(\#{\rm Bins} - 1)-1-3 = 20$ degrees of freedom
(the factor of $3$ which we subtract arises from estimating three
parameters, $\ga_\rs$, $\ga_\ra$ and $\gamma$; $\beta$ was not
estimated, as it was taken to be $-.5$).}. For the Weibulls
approximating the runs scored and allowed per game we used the bins
of \eqref{eq:listbaseballbins}: \be [-.5, .5] \ \cup \ [.5, 1.5] \
\cup \ \cdots \ \cup \ [7.5, 8.5]\ \cup\ [8.5,9.5]\ \cup\ [9.5,11.5]
\ \cup\ [11.5, \infty). \ee There are $20$ degrees of freedom for
these tests. For $20$ degrees of freedom the critical thresholds are
31.41 (at the $95\%$ level) and 37.57 (at the $99\%$ level).

We also tested the independence of the runs scored and runs allowed
per game (a crucial input for our model). As this test requires each
row and column to have at least one non-zero entry, here we broke
the runs scored and allowed into bins \be\label{eq:binsindepchitest}
[0, 1) \ \cup \ [1,2) \ \cup \ [2,3) \ \cup \cdots \ \cup \ [8, 9) \
\cup \ [9,10) \ \cup \ [10,11) \ \cup \ [11,\infty). \ee This gives
us an $r \times c$ contingency table (with $r=c=12$); however, as
the runs scored and allowed per game can never be equal, we actually
have an incomplete two-dimensional contingency table with $(12-1)^2
- 12 = 109$ degrees of freedom; see \cite{BF,SD}. This complication
is not present in the theoretical model, as if the runs scored and
allowed are drawn from continuous distributions (in this case,
Weibulls), there is zero probability of both values being equal.
This difficulty is due to the fact that the runs scored and allowed
in a game must be distinct integers. We describe the modified
$\chi^2$ test for an incomplete two-dimensional contingency table
with diagonal entries forced to be zero (these are called structural
or fixed zeros).

Let ${\rm Bin}(k)$ denote the $k$\textsuperscript{th} bin in
\eqref{eq:binsindepchitest}. For our $12\times 12$ incomplete
contingency table with these bins for both runs scored and allowed,
the entry $O_{r,c}$ corresponds to the observed number of games
where the team's runs scored is in ${\rm Bin}(r)$ and the runs
allowed are in ${\rm Bin}(c)$; note\footnote{The reason $O_{r,r}$
should equal zero is that a team cannot score and allow the same
number of runs in a game, as baseball does not allow ties (except
for an occasional All-star game). The first 11 bins each contain
exactly one score, so for $r \le 11$, $O_{r,r} = 0$. The final bin,
however, contains all scores from 11 to $\infty$, and thus it is
possible for the runs scored and allowed to be unequal and both in
this bin; however, the probability is so small here that we may
simply replace all runs scored or allowed exceeding 11 with 11. Of
the 14 teams, 7 have $O_{12,12} = 0$, 5 (teams 3, 4, 5, 9 and 12)
have $O_{12,12} = 1$, 1 (team 8) has $O_{12,12} = 2$ and 1 (team 7)
has $O_{12,12} = 3$.} $O_{r,r} = 0$ for all $r$. We use the
iterative fitting procedure given in the appendix to \cite{BF} to
obtain maximum likelihood estimators for the $E_{r,c}$, the expected
frequency of cell $(r,c)$ under the assumption that, given that the
runs scored and allowed are distinct, the runs scored and allowed
are independent. For $1 \le r,c \le 12$, let $E_{r,c}^{(0)} = 1$ if
$r \neq c$ and $0$ if $r=c$. Set \be X_{r,+} \ = \ \sum_{c=1}^{12}
O_{r,c}, \ \ \ \ X_{+,c} \ = \ \sum_{r=1}^{12} O_{r,c}. \ee Then \be
\threecase{E_{r,c}^{(\ell)} \ = \ }{E_{r,c}^{(\ell-1)} X_{r,+}\ /\
\sum_{c=1}^{12} E_{r,c}^{(\ell-1)}}{if $\ell$  is odd}{\ }{\
}{E_{r,c}^{(\ell-1)} X_{+,c}\ /\ \sum_{r=1}^{12}
E_{r,c}^{(\ell-1)}}{if $\ell$ is even,} \ee and \be E_{r,c} \ = \
\lim_{\ell \to \infty} E_{r,c}^{(\ell)}; \ee the iterations converge
very quickly in practice\footnote{If we had a complete
two-dimensional contingency table, then the iteration reduces to the
standard values, namely $E_{r,c} = \sum_{c'} O_{r,c'} \cdot
\sum_{r'} O_{r',c} \ / \ \#{\rm Games}$.}. Then \be \sum_{r=1}^{12}\
\sum_{c=1 \atop c \neq r}^{12}\ \frac{(O_{r,c} -
E_{r,c})^2}{E_{r,c}} \ee is approximately a $\chi^2$ distribution
with $(12-1)^2-12 = 109$ degrees of freedom. The corresponding
critical thresholds are 134.4 (at the $95\%$ level) and 146.3 (at
the $99\%$ level).


We summarize our results below; the first column is the $\chi^2$
tests for the goodness of fit from the best fit Weibulls, and the
second column is the $\chi^2$ tests for the independence of the runs
scored and runs allowed.

\begin{center} Results from the Method of Maximum Likelihood
\scalebox{.7}{\includegraphics{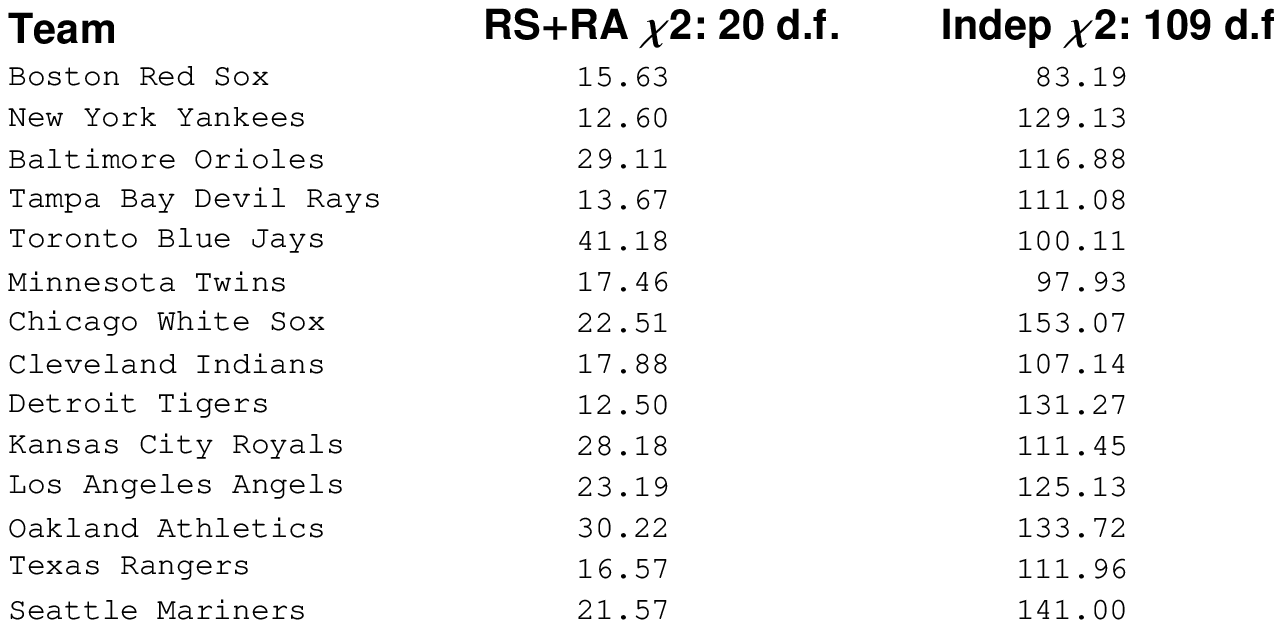}}
\end{center}

Except for the Weibulls for the runs scored and allowed for the
Toronto Blue Jays, and the independence of runs scored and runs
allowed for the Chicago White Sox\footnote{The Chicago White Sox had
the largest value of $O_{12,12}$ in the independence tests, namely
$3$. If we replace the last bin in \eqref{eq:binsindepchitest} with
two bins, $[11,12)$ and $[12,\infty)$, then $r=c=13$, $O_{12,12} =
0$ and $O_{13,13} = 1$. There are $(13-1)^2-13 = 131$ degrees of
freedom. The corresponding critical thresholds are 158.7 (at the
$95\%$ level) and 171.6 (at the $99\%$ level), and the observed
value of the $\chi^2$ statistic for the Chicago White Sox is
164.8.}, all test statistics are well below the $95\%$ critical
threshold (31.41 as there are 20 degrees of freedom). As we are
performing multiple comparisons, chance fluctuations should make
some differences appear significant (for example, if the null
hypothesis is true and 10 independent tests are performed, there is
about a $40\%$ chance of observing at least one statistically
significant difference at the $95\%$ confidence level). We must
therefore adjust the confidence levels. Using the common, albeit
conservative, Bonferroni\footnote{Using the Bonferroni adjustment
for multiple comparisons divides the significance level $\ga$ by the
number of comparisons, which in our case is $14$. Thus for the
Weibull tests with 20 degrees of freedom the adjusted critical
thresholds are 41.14 (at the $95\%$ level) and 46.38 (at the $99\%$
level); for the independence tests with 109 degrees of freedom the
adjusted critical thresholds are 152.9 (at the $95\%$ level) and
162.2 (at the $99\%$ level).} adjustment method for multiple
comparisons, at the $95\%$ confidence level we find significant fits
for all but the Toronto Blue Jays' runs scored and allowed and the
independence of runs scored and allowed for the Chicago White Sox;
however, both just barely miss at the $95\%$ confidence level (41.18
versus 41.14 for the Blue Jays, and 153.07 versus 152.9 for the
White Sox). Thus the data validates our assumption that, given that
runs scored and allowed cannot be equal, the runs scored and allowed
per game are statistically independent events, and that the
parameters from the method of maximum likelihood give good fits to
the observed distribution of scores. In Appendix
\ref{sec:plotsWeibulls} we provide plots comparing the observed
distribution of runs scored and allowed versus the best fit
predictions, where even a visual inspection shows the agreement
between our theory and the data.

Using the best fit parameters of the Weibulls, Lemma
\ref{lem:Weibullmeanvar} provides an estimate for the mean number of
runs scored and allowed per game. We are of course primarily
interested in estimating $\g$ and not the mean number of runs scored
or allowed per game, because these are of course known from the
season data; however, this provides an additional test to see how
well our theory agrees with the data.

As the number of games each team played is so large\footnote{All
teams played $162$ except for the Tampa Bay Devil Rays and the
Toronto Blue Jays, who had a game rained out and only played 161
games in 2004.}, we use a $z$-test to compare the observed versus
predicted means. The critical $z$-values are 1.96 (at the $95\%$
confidence level) and 2.575 (at the $99\%$ confidence level).

\begin{center} Results from the Method of Maximum Likelihood
\scalebox{.7}{\includegraphics{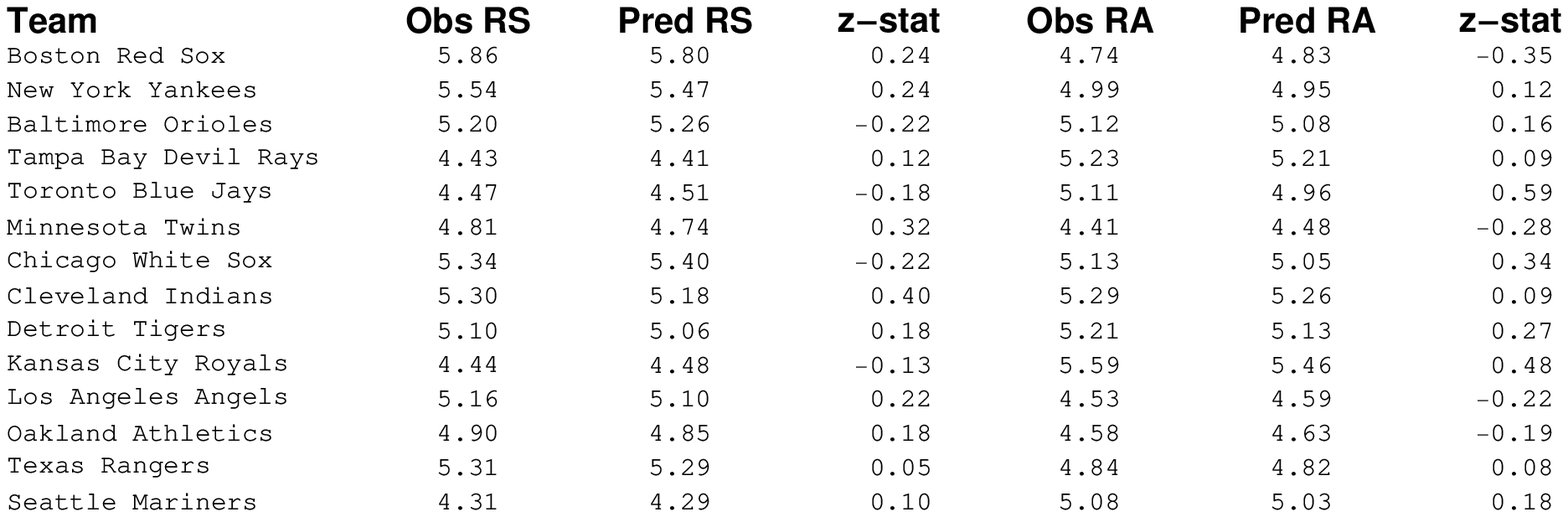}}
\end{center}

We note excellent agreement between all the predicted average runs
scored per game and the observed average runs scored per game, as
well as between all the predicted average runs allowed per game
and the observed average runs allowed per game. Performing a
Bonferroni adjustment for multiple comparisons gives critical
thresholds of 2.914 (at the $95\%$ level) and 3.384 (at the $99\%$
level). At the $95\%$ level (resp. $99\%$ level) all 14 teams have
significant fits.

As a final experiment, instead of finding the best fit Weibulls team
by team, we performed a similar analysis for each division in the
American League in 2004. For example, in the AL East there are 5
teams (the World Champion Boston Red Sox, the New York Yankees, the
Baltimore Orioles, the Tampa Bay Devil Rays and the Toronto Blue
Jays), and we found the least squares fit to the data with the 11
free parameters \be \ga_{\rs,{\rm BOS}},\ \ga_{\ra,{\rm BOS}}, \
\dots, \ \ga_{\rs,{\rm TOR}},\ \ga_{\ra,{\rm TOR}},\ \g. \ee The
five teams in the AL East (resp., the five teams of the AL Central
and the four teams of the AL West) give a best fit value of $\g$ of
$1.793$ (resp., 1.773 and 1.774), which again is very close to the
numerically observed best value of $\g$ of 1.82. Using the method of
maximum likelihood gives best fit values of $\gamma$ of 1.74 for the
AL East, $1.75$ for the AL Central and $1.73$ for the AL West.

\section{Conclusions and Future Work}\label{sec:concfutwork}

Bill James' Pythagorean Won-Loss Formula may be derived from very
simple and reasonable assumptions (namely, that the runs scored and
allowed per game are independent events drawn from Weibulls with the
same $\gb$ and $\g$). Using the method of least squares or the
method of maximum likelihood, we can find the best fit values of
these parameters from the observed game scores. Using the method of
maximum likelihood, for the 2004 baseball season for each team in
the American League the fits were always significant at the $95\%$
confidence level (except for the Toronto Blue Jays, which just
missed), the assumption that, given that the runs scored and allowed
in a game are distinct integers, the runs scored and allowed per
game are independent events was validated, and the best fit exponent
$\g$ was about $1.74$ with a standard deviation of $.06$, in
excellent agreement with the observation that $1.82$ is the best
exponent to use in the Pythagorean Formula (the method of least
squares gives a best fit value for $\g$ of $1.79$ with standard
deviation $.09$). Note that we obtain our value of the exponent $\g$
not by fitting the Pythagorean Formula to the observed won-loss
percentages of teams, but rather from an analysis of the
distribution of scores from individual baseball games. Assuming
teams behave similarly from year to year, there is now a theoretical
justification for using the Pythagorean Formula to predict team
performances in future seasons (with an exponent around $1.74$ to
$1.79$ and using the observed average runs scored and allowed).

An interesting future project would be to perform a more micro
analysis to incorporate lower order effects, though our simple model
is quite effective at fitting the data and predicting the best
exponent $\g$ (see for example \cite{Schell1,Schell2}, where such an
analysis is performed to determine the all-time best hitters and
sluggers). For example, one might break down runs scored and allowed
per inning. If a team has a large lead it often pulls its good
hitters to give them a rest, as well as bringing in weaker pitchers
to give them some work; conversely, in late innings in close games
managers often have pitch-runners for slow good hitters who get on
base (to get a run now with a potential cost of runs later through
the loss of the better hitter from the lineup), and star relievers
(when available) are brought in to maintain the lead. Further there
are slight differences because of inter-league play. For example,
the American League teams lose their DH for games in National League
parks, and thus we expect the run production to differ from that in
American League parks. Further, using the analysis in
\cite{Schell1,Schell2} one can incorporate ballpark effects (some
ballparks favor pitchers while others favor hitters). Such an
analysis might lead to new statistics of adjusted runs scored and
allowed per game. Additionally, teams out of the playoff race often
play their last few games differently than when they are still in
contention, and perhaps those games should be removed from the
analysis.

One can also further examine the independence of runs scored and
allowed. As baseball games cannot end in a tie, runs scored and
allowed are never equal in a game; however, they can be equal after
9 innings. One avenue for research is to classify extra-inning games
as ties (while recording which team eventually won). Also, if the
home team is leading after the top of the ninth then it does not
bat, and this will effect its run production. See \cite{Ci} for an
analysis of some of these issues.

Finally, it would be fascinating to see if this (or a similar)
model is applicable to other sports with long seasons. While
football has a relatively short season of 16 games, basketball and
hockey have 82 games a season. The scores in basketball are more
spread out than hockey, which is more compact than baseball; it
would be interesting to see what affect these have on the analysis
and whether or not the fits are as good as baseball.

\section*{Acknowledgements}

I would like to thank Russell Mann and Steven Johnson for
introducing me to the Pythagorean Formula in baseball, Kevin
Dayaratna for inputting much of the baseball data, Jeff Miller for
writing a script to read in baseball data from the web to the
analysis programs, Gerry Myerson for catching some typos, Eric T.
Bradlow for helpful comments on an earlier draft, Ray Ciccolella for
discussions on the independence of runs scored and allowed, and
Stephen D. Miller for suggesting the National League exercise.


\appendix

\section{Plots of Best Fit Weibulls}\label{sec:plotsWeibulls}

Below we plot the best fit Weibulls against the observed histograms
of runs scored and allowed. We use the bins of
\eqref{eq:listbaseballbins}.


\begin{center}
\scalebox{.6}{\includegraphics{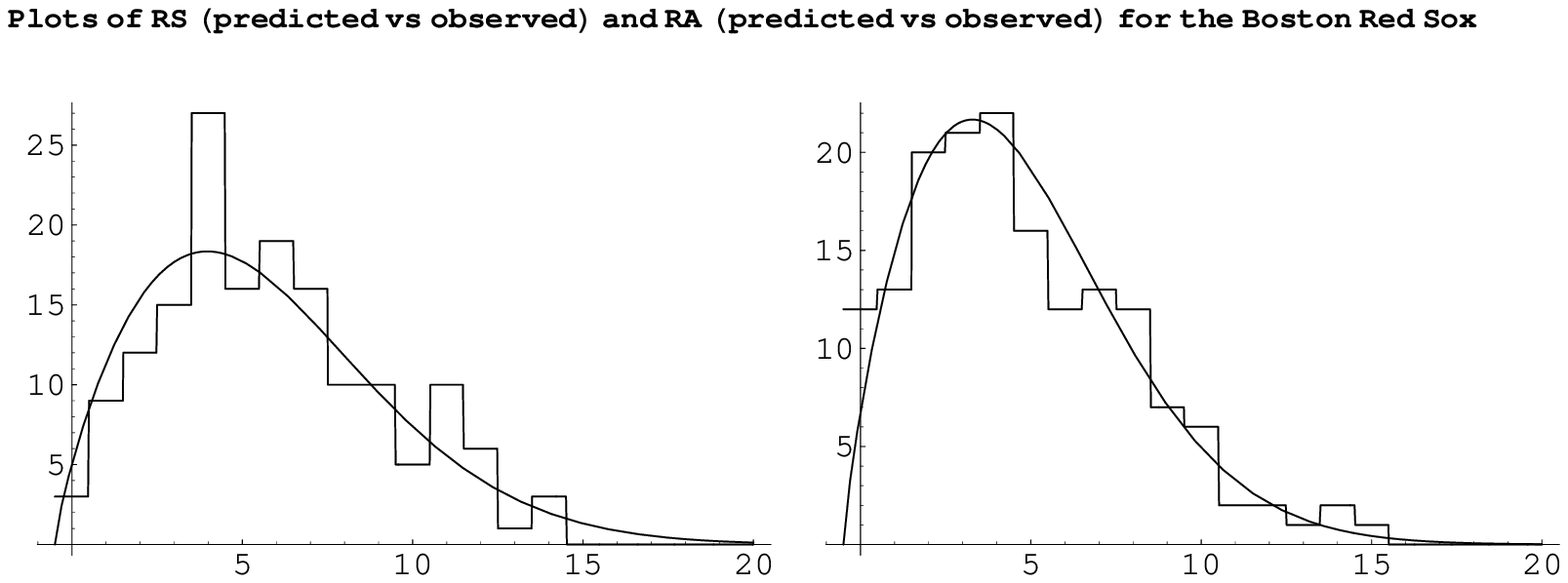}}
\end{center}

\bigskip

\begin{center}
\scalebox{.6}{\includegraphics{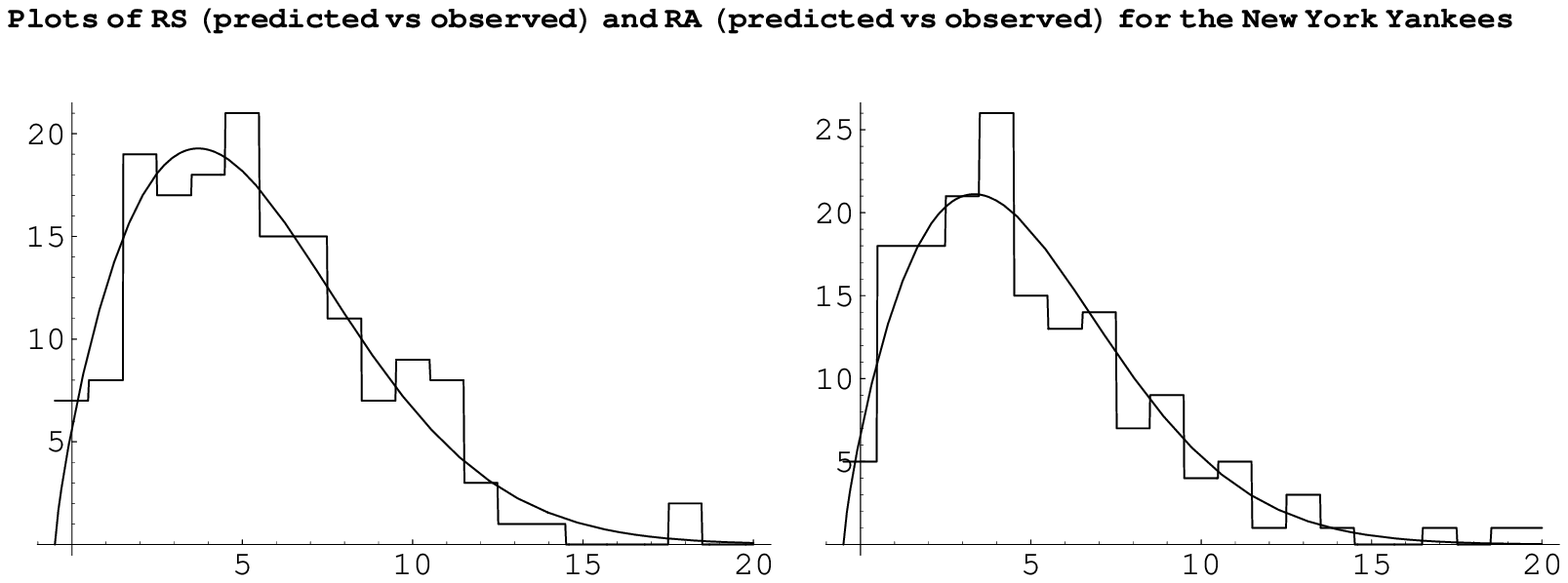}}
\end{center}

\bigskip

\begin{center}
\scalebox{.6}{\includegraphics{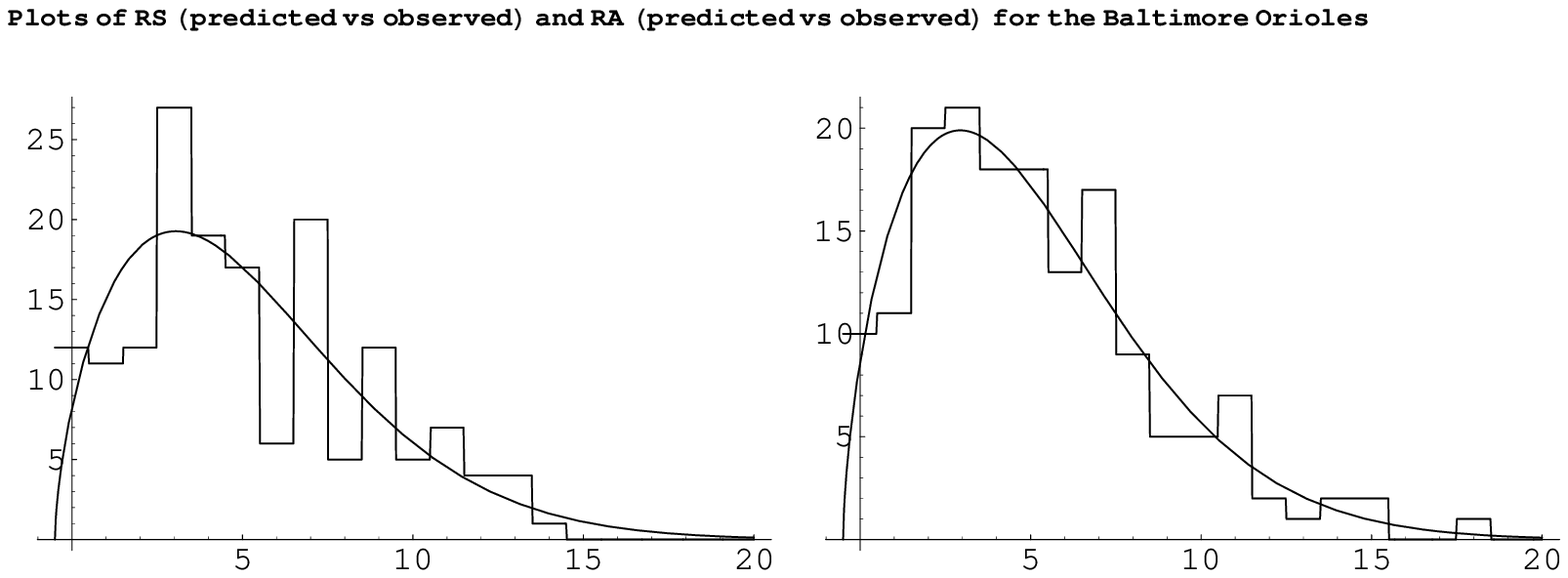}}
\end{center}

\bigskip

\begin{center}
\scalebox{.6}{\includegraphics{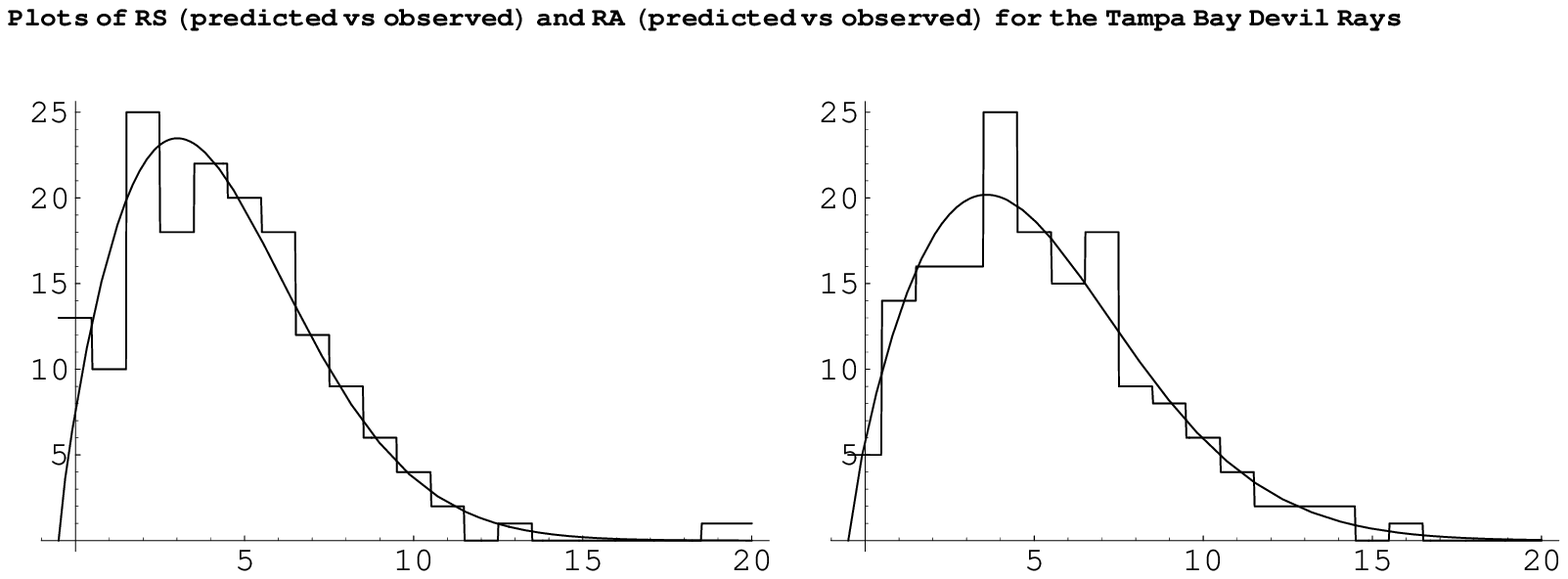}}
\end{center}

\bigskip

\begin{center}
\scalebox{.6}{\includegraphics{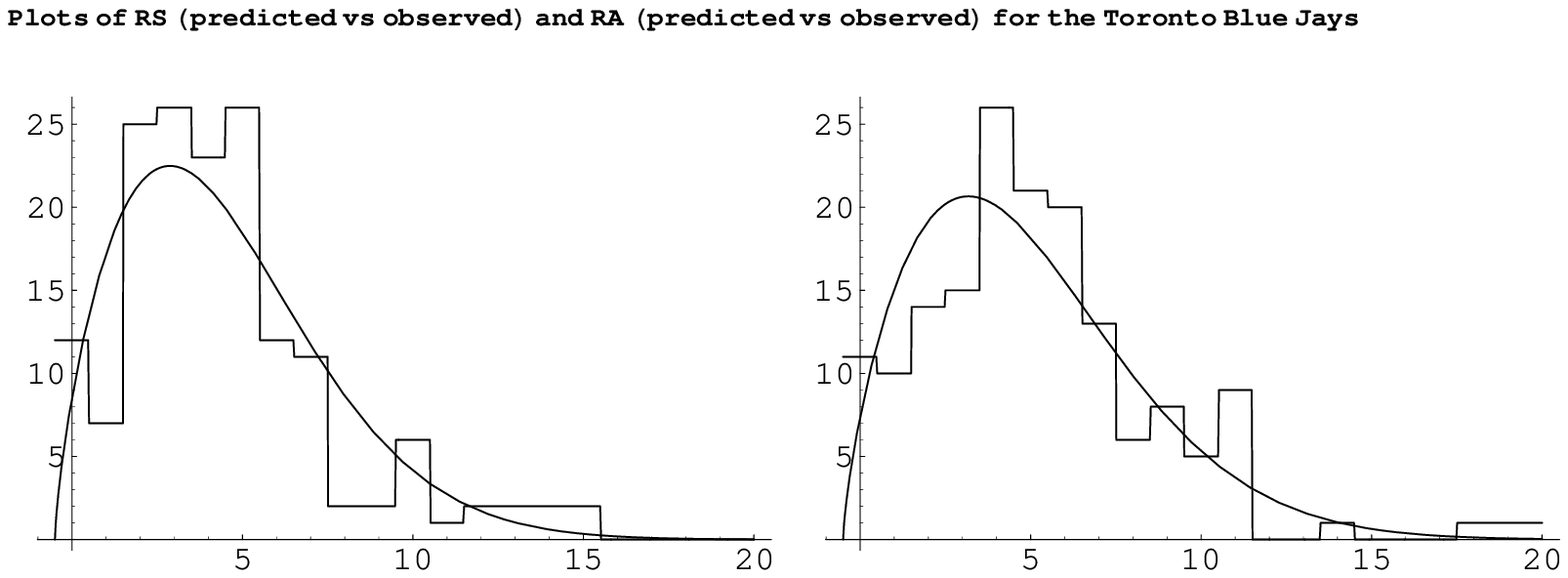}}
\end{center}

\bigskip


\begin{center}
\scalebox{.6}{\includegraphics{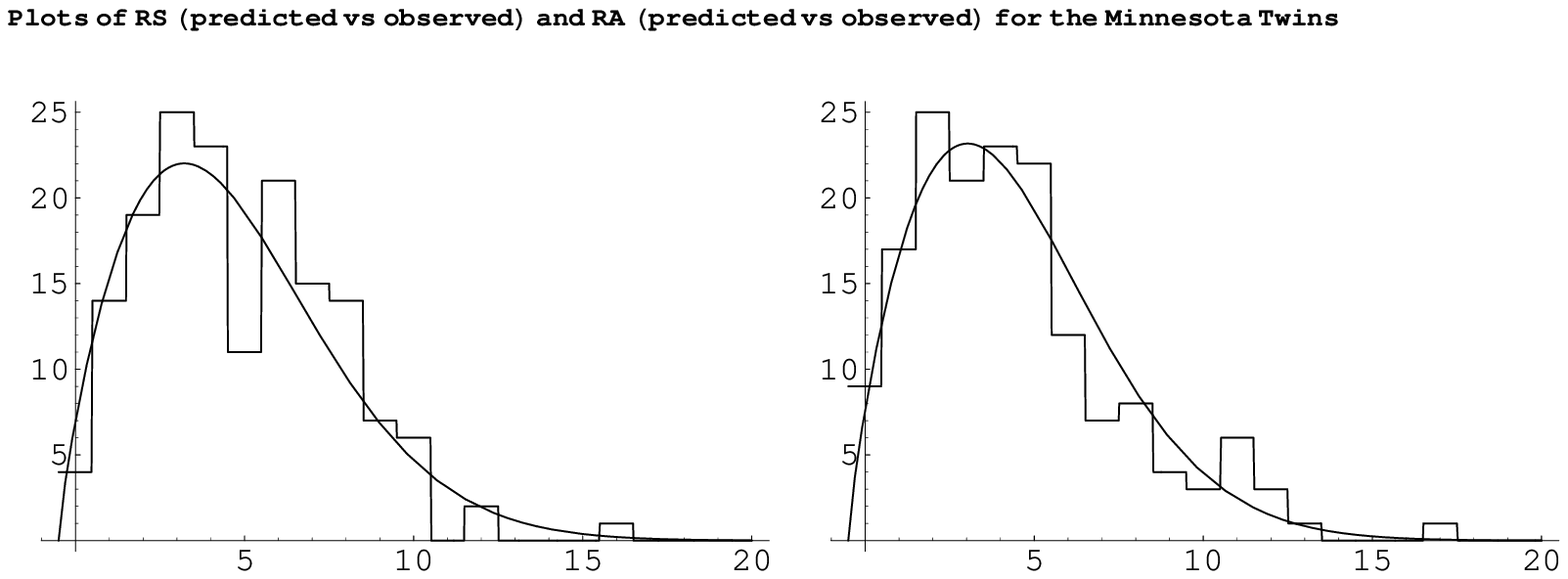}}
\end{center}

\bigskip

\begin{center}
\scalebox{.6}{\includegraphics{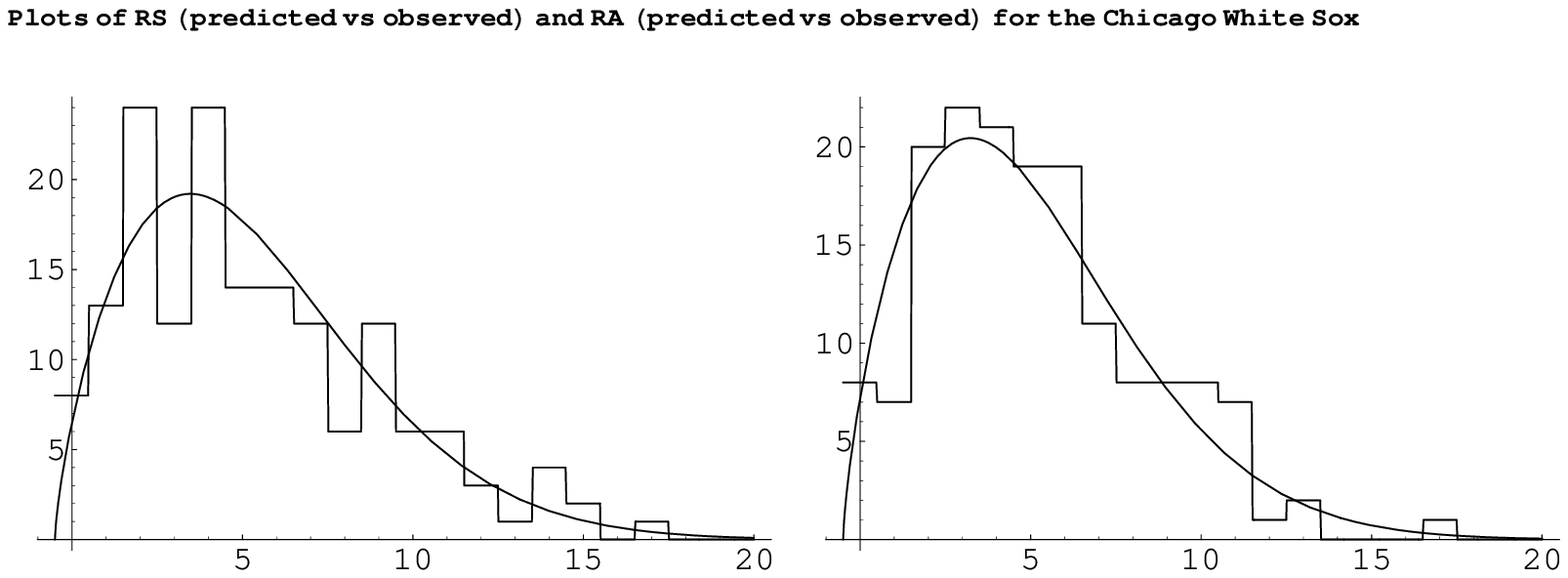}}
\end{center}

\bigskip

\begin{center}
\scalebox{.6}{\includegraphics{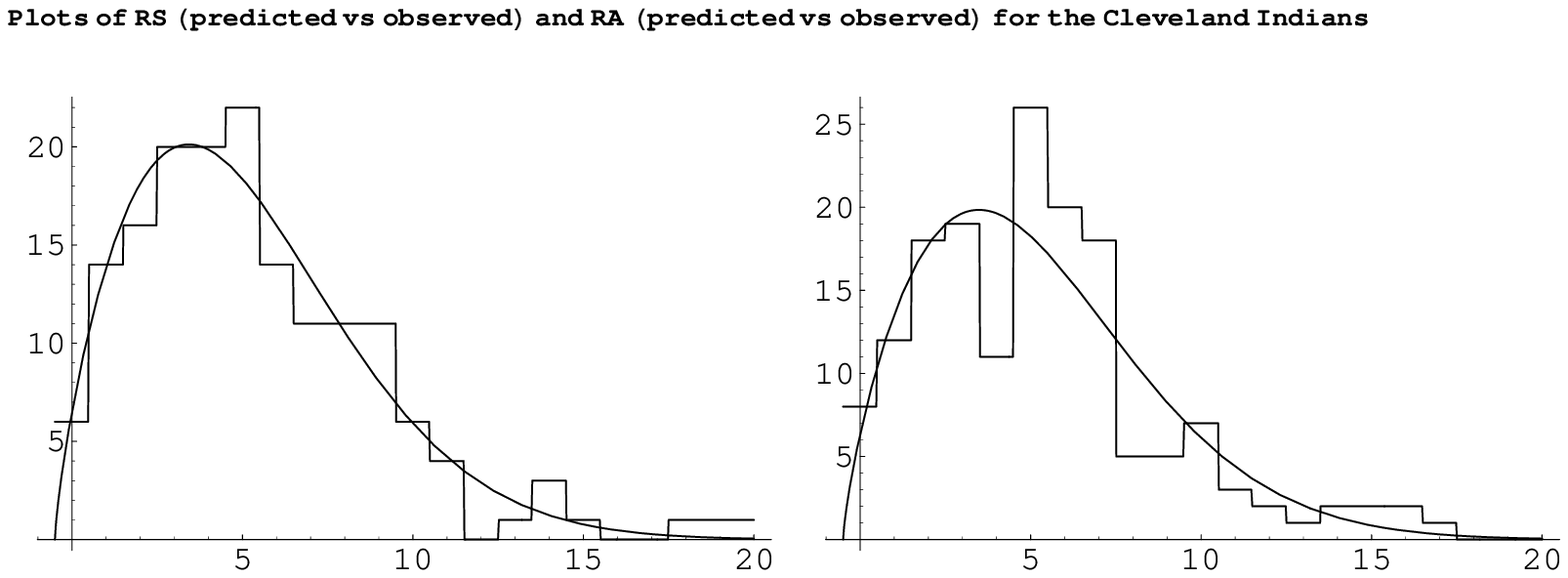}}
\end{center}

\bigskip

\begin{center}
\scalebox{.6}{\includegraphics{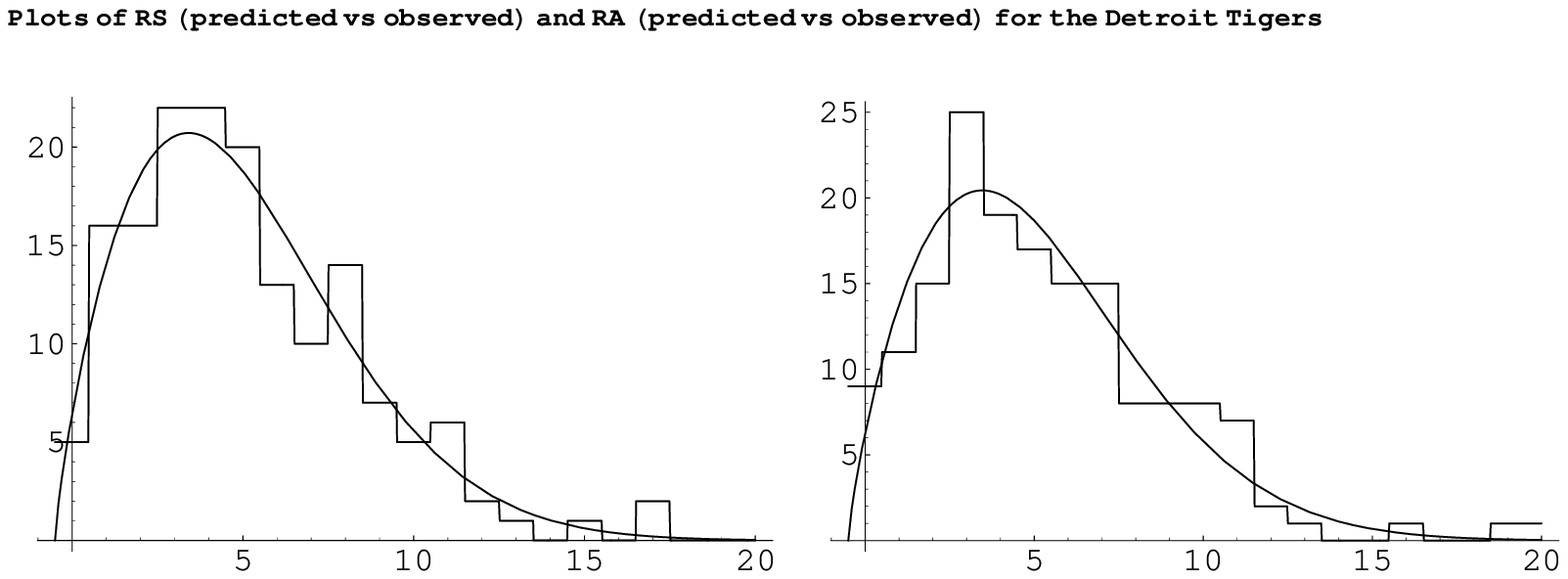}}
\end{center}

\bigskip

\begin{center}
\scalebox{.6}{\includegraphics{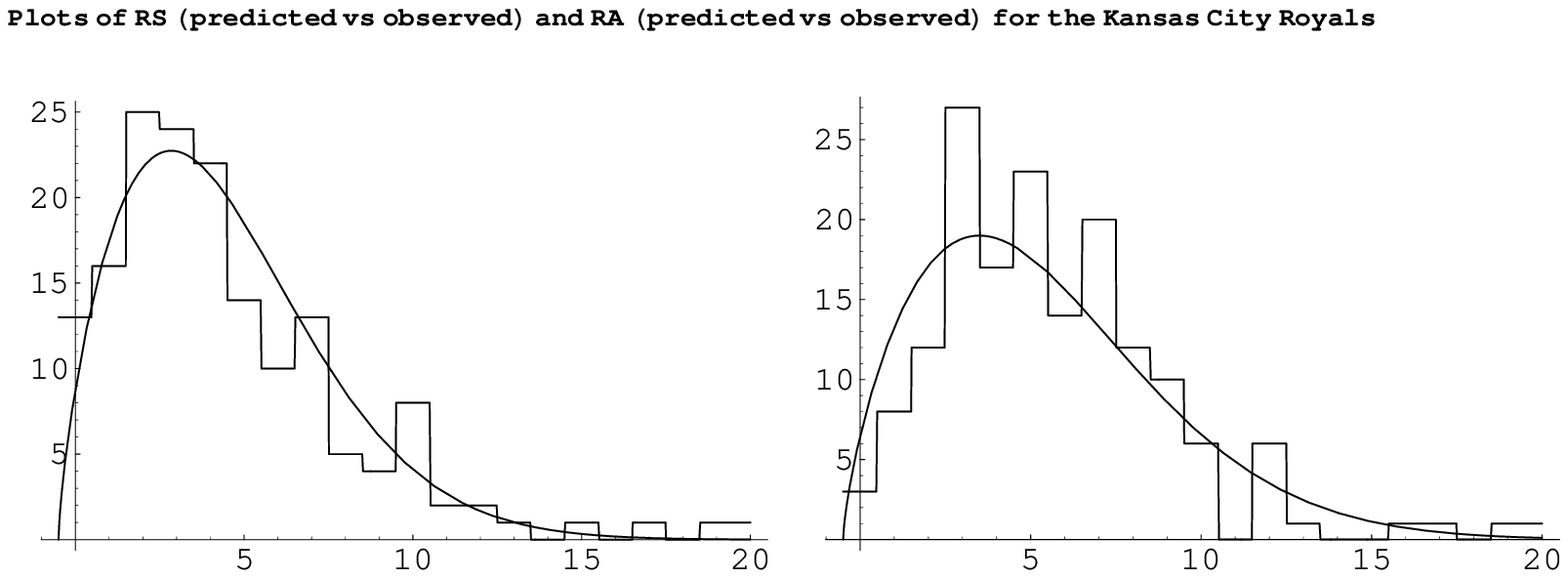}}
\end{center}

\bigskip


\end{document}